\documentstyle{article}

\def\cirk{\,{\raisebox{.3ex}{\tiny $\circ$}}\,}
\newcommand{\mj}{\mbox{\bf 1}}

\begin{document}

\title{Coherence and Confluence}
\author{{\sc Kosta Do\v sen} and {\sc Zoran Petri\' c}\\[0.5cm]
Mathematical Institute, SANU \\
Knez Mihailova 35, p.f. 367 \\
11001 Belgrade, Serbia \\
email: \{kosta, zpetric\}@mi.sanu.ac.yu}
\date{}
\maketitle

\begin{abstract}
\noindent Proofs of coherence in category theory, starting from
Mac Lane's original proof of coherence for monoidal categories,
are sometimes based on confluence techniques analogous to what one
finds in the lambda calculus, or in term-rewriting systems in
general. This applies to coherence results that assert that a
category is a preorder, i.e. that ``all diagrams commute''. This
note is about this analogy, paying particular attention to cases
where the category for which coherence is proved is not a
groupoid.
\end{abstract}

\vspace{.3cm}

{\small 18A15, 18D10, 03B40, 03D03}

\vspace{.5ex}

\noindent {\it Keywords$\,$}: {\small coherence, categories,
confluence, Church-Rosser property, term-rewriting systems,
equational theories}

\baselineskip=1.075\baselineskip

\section{Introduction}

This note is about a connection between the categorial notion of
coherence and the notion of confluence found in term-rewriting
systems. By coherence we understand the following:

\begin{quotation}

\noindent {\it Coherence is a completeness result for an
axiomatization of a brand of category, usually with respect to a
particular category as a model.}

\end{quotation}

In cases when one expects from coherence to decide whether two
terms stand for the same arrow, the model category should be
manageable in the sense that there is a decision procedure,
preferably elementary, for equality of arrows in it. By varying
the model category, we can cover with the notion above the results
of Mac Lane and Kelly concerning coherence of monoidal, symmetric
monoidal and symmetric monoidal closed categories (see
\cite{ML63}, \cite{ML71} and \cite{KML71}), as well as many other
coherence results (see \cite{DP04} and \cite{DP05}).

This notion of coherence is made more precise by taking in the
particular brand of category that interests us a category $\cal K$
freely generated by a set of objects (this set may be understood
as a discrete category). This free category will always exist if
our axiomatization is purely equational. Then coherence amounts to
showing the following:

\begin{quotation}

\noindent {\it There is a faithful functor $G$ from the free
category $\cal K$ to a particular model category $\cal M$.}

\end{quotation}

\noindent In logical terms, the existence of the functor $G$ from
$\cal K$ to $\cal M$ is soundness, and the faithfulness of $G$ is
completeness proper.

Proofs of coherence in category theory, starting from Mac Lane's
original proof of coherence for monoidal categories of
\cite{ML63}, are sometimes based on confluence techniques
analogous to what one finds in the lambda calculus, or in
term-rewriting systems in general. This applies to coherence
results that assert that a category is a preorder, i.e. that ``all
diagrams commute''. (A preordering relation is a reflexive and
transitive relation; a category that is a preorder is a
preordering relation on the set of its objects.) To make such
coherence results accord with the notion of coherence above, in
many cases one can take that the image of $\cal K$ in $\cal M$ is
a discrete category. In this note we will make some comments on
the analogy between proofs of coherence and proofs of confluence,
paying particular attention to cases where the category for which
coherence is proved is not a groupoid.

\section{Coherence and proof theory}
If one envisages a deductive system as a graph whose nodes are
formulae:

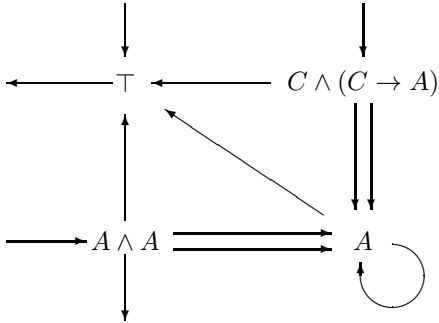
\begin{center}
\begin{picture}(150,120)
\put(45,30){\makebox(0,0){$A\wedge A$}}
\put(45,90){\makebox(0,0){$\top$}}
\put(135,30){\makebox(0,0){$A$}}
\put(135,90){\makebox(0,0){$C\wedge (C\rightarrow A)$}}

\put(0,30){\vector(1,0){30}}
\put(63,27){\vector(1,0){60}}
\put(63,33){\vector(1,0){60}}
\put(45,25){\vector(0,-1){25}}
\put(45,38){\vector(0,1){40}}
\put(132,82){\vector(0,-1){40}}
\put(138,82){\vector(0,-1){40}}
\put(120,40){\vector(-3,2){60}}
\put(40,90){\vector(-1,0){40}}
\put(100,90){\vector(-1,0){45}}
\put(45,120){\vector(0,-1){20}}
\put(135,120){\vector(0,-1){20}}

\put(146,17){\oval(24,24)[r]}
\put(146,17){\oval(24,24)[bl]}
\put(134,17){\vector(0,1){5}}

\end{picture}
\end{center}

\noindent and whose arrows are derivations from the sources
understood as premises to the targets understood as conclusions,
then equality of derivations usually transforms this deductive
system into a category of a particular brand. This category has a
structure induced by the connectives of the deductive system.
Although equality of derivation is dictated by logical concerns,
usually the categories we end up with are of a kind that
categorists have already introduced for their own reason. The
prime example here is given by the deductive system for the
conjunction-implication fragment of intuitionistic propositional
logic. After derivations in this deductive system are equated
according to ideas about normalization of derivations that stem
from Gentzen, one obtains the cartesian closed category $\cal K$
freely generated by a set of propositional variables.

Equality of proofs in intuitionistic logic has not led up to now
to a coherence result---a coherence theorem is not forthcoming for
cartesian closed categories. If we take that the model category
$\cal M$ is a category whose arrows are graphs like the graphs of
\cite{KML71}, then we do not have a faithful functor $G$ from the
free cartesian closed category $\cal K$ to $\cal M$.

If $\eta_{p,q}$ is the canonical arrow from $q$ to $p\rightarrow
(p\times q)$, where $A\rightarrow B$ stands for $B^A$, and $w_A$
is the diagonal arrow from $A$ to $A\times A$, then
$G(w_{p\rightarrow (p\times q)}\cirk\eta_{p,q})$:

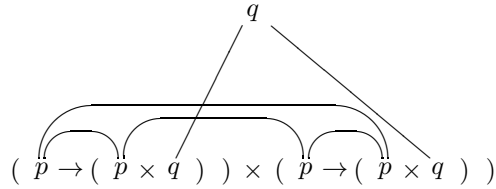
\begin{center}
\begin{picture}(180,80)
\put(10,10){\makebox(0,0){$($}} \put(20,10){\makebox(0,0){$p$}}
\put(31,10){\makebox(0,0){$\rightarrow$}}
\put(40,10){\makebox(0,0){$($}} \put(50,10){\makebox(0,0){$p$}}
\put(60,10){\makebox(0,0){$\times$}}
\put(70,10){\makebox(0,0){$q$}} \put(80,10){\makebox(0,0){$)$}}
\put(90,10){\makebox(0,0){$)$}}
\put(100,10){\makebox(0,0){$\times$}}
\put(110,10){\makebox(0,0){$($}} \put(120,10){\makebox(0,0){$p$}}
\put(131,10){\makebox(0,0){$\rightarrow$}}
\put(140,10){\makebox(0,0){$($}} \put(150,10){\makebox(0,0){$p$}}
\put(160,10){\makebox(0,0){$\times$}}
\put(170,10){\makebox(0,0){$q$}} \put(180,10){\makebox(0,0){$)$}}
\put(190,10){\makebox(0,0){$)$}}

\put(100,70){\makebox(0,0){$q$}}

\put(35,15){\oval(28,20)[t]}

\put(85,15){\oval(132,40)[t]}

\put(135,15){\oval(28,20)[t]}

\put(85,15){\oval(68,30)[t]}

\put(71,15){\line(1,2){25}}

\put(167,15){\line(-6,5){60}}

\end{picture}
\end{center}

\noindent which is obtained from

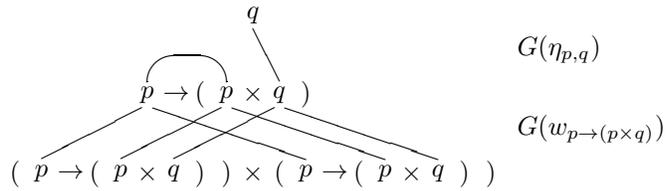
\begin{center}
\begin{picture}(180,80)
\put(10,10){\makebox(0,0){$($}} \put(20,10){\makebox(0,0){$p$}}
\put(31,10){\makebox(0,0){$\rightarrow$}}
\put(40,10){\makebox(0,0){$($}} \put(50,10){\makebox(0,0){$p$}}
\put(60,10){\makebox(0,0){$\times$}}
\put(70,10){\makebox(0,0){$q$}} \put(80,10){\makebox(0,0){$)$}}
\put(90,10){\makebox(0,0){$)$}}
\put(100,10){\makebox(0,0){$\times$}}
\put(110,10){\makebox(0,0){$($}} \put(120,10){\makebox(0,0){$p$}}
\put(131,10){\makebox(0,0){$\rightarrow$}}
\put(140,10){\makebox(0,0){$($}} \put(150,10){\makebox(0,0){$p$}}
\put(160,10){\makebox(0,0){$\times$}}
\put(170,10){\makebox(0,0){$q$}} \put(180,10){\makebox(0,0){$)$}}
\put(190,10){\makebox(0,0){$)$}}

\put(60,40){\makebox(0,0){$p$}}
\put(71,40){\makebox(0,0){$\rightarrow$}}
\put(80,40){\makebox(0,0){$($}} \put(90,40){\makebox(0,0){$p$}}
\put(100,40){\makebox(0,0){$\times$}}
\put(110,40){\makebox(0,0){$q$}} \put(120,40){\makebox(0,0){$)$}}

\put(100,70){\makebox(0,0){$q$}}

\put(100,65){\line(1,-2){10}} \put(75,45){\oval(30,20)[t]}

\put(58,35){\line(-2,-1){38}} \put(88,35){\line(-2,-1){38}}
\put(108,35){\line(-2,-1){38}} \put(62,35){\line(3,-1){58}}
\put(92,35){\line(3,-1){58}} \put(112,35){\line(3,-1){58}}

\put(200,25){$G(w_{p\rightarrow (p\times q)})$}

\put(200,55){$G(\eta_{p,q})$}

\end{picture}
\end{center}

\noindent is different from $G((\eta_{p,q}\times\eta_{p,q})\cirk
w_q)$:

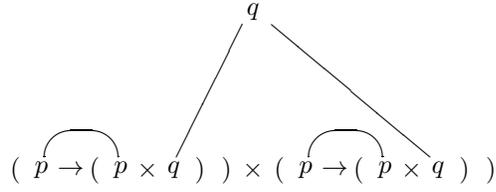
\begin{center}
\begin{picture}(180,80)
\put(10,10){\makebox(0,0){$($}} \put(20,10){\makebox(0,0){$p$}}
\put(31,10){\makebox(0,0){$\rightarrow$}}
\put(40,10){\makebox(0,0){$($}} \put(50,10){\makebox(0,0){$p$}}
\put(60,10){\makebox(0,0){$\times$}}
\put(70,10){\makebox(0,0){$q$}} \put(80,10){\makebox(0,0){$)$}}
\put(90,10){\makebox(0,0){$)$}}
\put(100,10){\makebox(0,0){$\times$}}
\put(110,10){\makebox(0,0){$($}} \put(120,10){\makebox(0,0){$p$}}
\put(131,10){\makebox(0,0){$\rightarrow$}}
\put(140,10){\makebox(0,0){$($}} \put(150,10){\makebox(0,0){$p$}}
\put(160,10){\makebox(0,0){$\times$}}
\put(170,10){\makebox(0,0){$q$}} \put(180,10){\makebox(0,0){$)$}}
\put(190,10){\makebox(0,0){$)$}}

\put(100,70){\makebox(0,0){$q$}}

\put(35,15){\oval(28,20)[t]}

\put(135,15){\oval(28,20)[t]}

\put(71,15){\line(1,2){25}}

\put(167,15){\line(-6,5){60}}

\end{picture}
\end{center}

\noindent which is obtained from

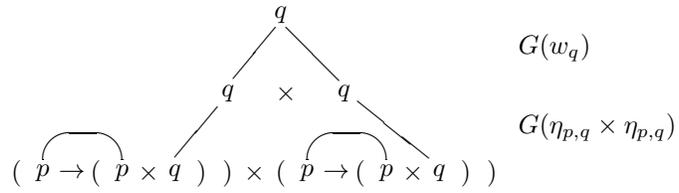
\begin{center}
\begin{picture}(400,80)(125,0)

\put(210,10){\makebox(0,0){$($}} \put(220,10){\makebox(0,0){$p$}}
\put(231,10){\makebox(0,0){$\rightarrow$}}
\put(240,10){\makebox(0,0){$($}} \put(250,10){\makebox(0,0){$p$}}
\put(260,10){\makebox(0,0){$\times$}}
\put(270,10){\makebox(0,0){$q$}} \put(280,10){\makebox(0,0){$)$}}
\put(290,10){\makebox(0,0){$)$}}
\put(300,10){\makebox(0,0){$\times$}}
\put(310,10){\makebox(0,0){$($}} \put(320,10){\makebox(0,0){$p$}}
\put(331,10){\makebox(0,0){$\rightarrow$}}
\put(340,10){\makebox(0,0){$($}} \put(350,10){\makebox(0,0){$p$}}
\put(360,10){\makebox(0,0){$\times$}}
\put(370,10){\makebox(0,0){$q$}} \put(380,10){\makebox(0,0){$)$}}
\put(390,10){\makebox(0,0){$)$}}

\put(290,40){\makebox(0,0){$q$}}
\put(312,40){\makebox(0,0){$\times$}}
\put(334,40){\makebox(0,0){$q$}}

\put(310,70){\makebox(0,0){$q$}}

\put(308,65){\line(-5,-6){16}} \put(312,65){\line(1,-1){20}}

\put(235,15){\oval(30,20)[t]} \put(335,15){\oval(30,20)[t]}

\put(286,35){\line(-5,-6){16}} \put(339,37){\line(5,-4){27}}

\put(400,25){$G(\eta_{p,q}\times\eta_{p,q})$}

\put(400,55){$G(w_q)$}

\end{picture}
\end{center}

\noindent So, if $w$ is a natural transformation, then $G$ is not
a functor.

Dually, if $\varepsilon_{p,q}$ is the canonical arrow from
$p\times(p\rightarrow q)$ to $q$, and $k^1_{A,B}$ is the first
projection from $A\times B$ to $A$, then
$G(k^1_{r,q}\cirk(\mj_r\times\varepsilon_{p,q}))$:

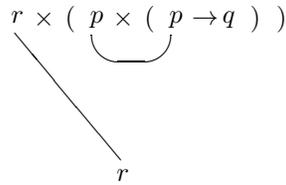
\begin{center}
\begin{picture}(120,80)

\put(50,11){\makebox(0,0){$r$}}

\put(10,71){\makebox(0,0){$r$}}

\put(20,70){\makebox(0,0){$\times$}}

\put(30,70){\makebox(0,0){$($}}

\put(40,70){\makebox(0,0){$p$}}

\put(50,70){\makebox(0,0){$\times$}}

\put(60,70){\makebox(0,0){$($}}

\put(70,70){\makebox(0,0){$p$}}

\put(81,70){\makebox(0,0){$\rightarrow$}}

\put(90,70){\makebox(0,0){$q$}}

\put(100,70){\makebox(0,0){$)$}}

\put(110,70){\makebox(0,0){$)$}}

\put(49,17){\line(-5,6){40}}

\put(53,64){\oval(30,20)[b]}

\end{picture}
\end{center}

\noindent which is obtained from

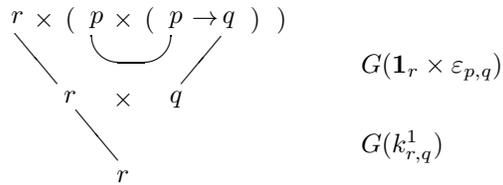
\begin{center}
\begin{picture}(120,80)

\put(50,11){\makebox(0,0){$r$}}

\put(30,41){\makebox(0,0){$r$}}

\put(50,40){\makebox(0,0){$\times$}}

\put(70,40){\makebox(0,0){$q$}}

\put(10,71){\makebox(0,0){$r$}}

\put(20,70){\makebox(0,0){$\times$}}

\put(30,70){\makebox(0,0){$($}}

\put(40,70){\makebox(0,0){$p$}}

\put(50,70){\makebox(0,0){$\times$}}

\put(60,70){\makebox(0,0){$($}}

\put(70,70){\makebox(0,0){$p$}}

\put(81,70){\makebox(0,0){$\rightarrow$}}

\put(90,70){\makebox(0,0){$q$}}

\put(100,70){\makebox(0,0){$)$}}

\put(110,70){\makebox(0,0){$)$}}

\put(48,17){\line(-5,6){16}}

\put(25,46){\line(-5,6){16}}

\put(72,46){\line(5,6){15}}

\put(53,64){\oval(30,20)[b]}

\put(140,20){$G(k^1_{r,q})$}

\put(140,50){$G(\mj_r\times\varepsilon_{p,q})$}

\end{picture}
\end{center}

\noindent is different from $G(k^1_{r,p\times(p\rightarrow q)})$:

\begin{center}
\begin{picture}(120,80)

\put(50,11){\makebox(0,0){$r$}}

\put(10,71){\makebox(0,0){$r$}}

\put(20,70){\makebox(0,0){$\times$}}

\put(30,70){\makebox(0,0){$($}}

\put(40,70){\makebox(0,0){$p$}}

\put(50,70){\makebox(0,0){$\times$}}

\put(60,70){\makebox(0,0){$($}}

\put(70,70){\makebox(0,0){$p$}}

\put(81,70){\makebox(0,0){$\rightarrow$}}

\put(90,70){\makebox(0,0){$q$}}

\put(100,70){\makebox(0,0){$)$}}

\put(110,70){\makebox(0,0){$)$}}

\put(49,17){\line(-5,6){40}}

\end{picture}
\end{center}

\noindent So, if $k^1$ is a natural transformation, then $G$ is
not a functor. The faithfulness of $G$ fails because of a
counterexample in \cite{S75}. This does not exclude that with a
more sophisticated model category $\cal M$ we might still be able
to obtain coherence for cartesian closed categories (for an
attempt along these lines see \cite{PD}).

Equality of proofs in classical logic may, however, lead to
coherence with respect to model categories that catch up to a
point the idea of \emph{generality} of proofs. Such is in
particular the category \emph{Rel}, whose arrows are relations
between finite ordinals, i.e.\ relations between occurrences of
the same propositional letters in the premises and conclusions.
The idea that generality of proofs may serve as a criterion for
identity of proofs stems from Lambek's pioneering papers in
categorial proof theory of the late 1960s (see \cite{LS86} for
references). This criterion says, roughly, that two derivations
represent the same proof when their generalizations with respect
to diversification of variables (without changing the rules of
inference) produce derivations with the same source and target, up
to a renaming of variables.

It is shown in \cite{DP04} that coherence with respect to the
model category \emph{Rel} could justify plausibly equality of
derivations in various systems of propositional logic, including
classical propositional logic. The goal of that book was to
explore the limits of coherence with respect to the model category
\emph{Rel}. This does not exclude that other coherence results may
involve other model categories, and, in particular, with a model
category different from \emph{Rel}, classical propositional logic
may induce a different notion of Boolean category than the one
introduced in Chapter 14 of \cite{DP04}. That notion of Boolean
category was not motivated \emph{a priori}, but was dictated by
coherence with respect to \emph{Rel}. The definition of that
notion was however not given via coherence, but via an equational
axiomatization. We take such definitions as being proper axiomatic
definitions.

We could easily define nonaxiomatically a notion of Boolean
category with respect to graphs of the Kelly-Mac Lane kind (see
\cite{KML71}). Equality of graphs would dictate what arrows are
equal. In this notion, conjunction would not be a product, because
the diagonal arrows and the projections would not make natural
transformations (see above), and, analogously, disjunction would
not be a coproduct (cf.\ \cite{DP04}, Section 14.3.) The resulting
notion of Boolean category would not be trivial---the freely
generated categories of that kind would not be preorders---, but
its nonaxiomatic definition would be trivial. There might exist a
nontrivial equational axiomatic definition of this notion. Finding
such a definition is an open problem.

We are looking for nontrivial axiomatic definitions because such
definitions give information about the combinatorial building
blocks of our notions, as Reidemeister moves give information
about the combinatorial building blocks of knot equivalence (see
\cite{BZ85}, Chapter 1). Our axiomatic equational definition of
Boolean category in \cite{DP04} is of the nontrivial,
combinatorially informative, kind. Coherence of these Boolean
categories with respect to \emph{Rel} is a theorem, whose proof in
\cite{DP04} requires considerable effort.

Another analogous example is provided by the notion of monoidal
category, which was introduced in a not entirely axiomatic way,
via coherence, by B\' enabou in \cite{Ben63}, and in the axiomatic
way, such as we favour, by Mac Lane in \cite{ML63}. For B\'
enabou, coherence is built into the definition, and for Mac Lane
it is a theorem. One could analogously define the theorems of
classical propositional logic as being the tautologies (this is
done, for example, in \cite{CK73}, Sections 1.2-3), in which case
completeness would not be a theorem, but would be built into the
definition.

\section{All diagrams commute}
The simplest case of coherence is when it asserts that ``all
diagrams commute'', which means that the free category $\cal K$ is
a preorder, i.e.\ a preordering relation on its objects. In this
case, some techniques used for proving coherence are related to
those developed in connection with term-rewriting systems (cf.\
\cite{Huet80} and \cite{JohnM92}). The difference is that with
coherence we are not interested in proving that starting from an
object all paths, i.e.\ all sequences, of arrows (reductions)
obtained by composing terminate in the same normal form. (This may
obtain sometimes, but is not essential.) Instead, we are
interested in proving that the equality of such paths follows from
some basic equations assumed for arrows. So, the level of our
interest is not the same. (This is why we need not go so high as
\cite{JohnM92} in the $n$-categorial hierarchy.)

Reductions here differ also from reductions in the lambda
calculus, where the lambda terms, which correspond to our arrows,
are reduced. We do not reduce arrows, but their types.

If all the arrows in question are isomorphisms, then proving that
all paths of arrows from the same source to the same target:

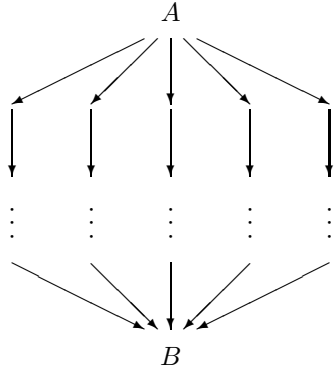
\begin{center}
\begin{picture}(140,140)
\put(70,5){\makebox(0,0){$B$}}
\put(70,135){\makebox(0,0){$A$}}

\put(10,40){\vector(2,-1){50}}
\put(40,40){\vector(1,-1){25}}
\put(70,40){\vector(0,-1){25}}
\put(100,40){\vector(-1,-1){25}}
\put(130,40){\vector(-2,-1){50}}

\put(10,50){\makebox(0,0){$\cdot$}}
\put(10,55){\makebox(0,0){$\cdot$}}
\put(10,60){\makebox(0,0){$\cdot$}}

\put(40,50){\makebox(0,0){$\cdot$}}
\put(40,55){\makebox(0,0){$\cdot$}}
\put(40,60){\makebox(0,0){$\cdot$}}

\put(70,50){\makebox(0,0){$\cdot$}}
\put(70,55){\makebox(0,0){$\cdot$}}
\put(70,60){\makebox(0,0){$\cdot$}}

\put(100,50){\makebox(0,0){$\cdot$}}
\put(100,55){\makebox(0,0){$\cdot$}}
\put(100,60){\makebox(0,0){$\cdot$}}

\put(130,50){\makebox(0,0){$\cdot$}}
\put(130,55){\makebox(0,0){$\cdot$}}
\put(130,60){\makebox(0,0){$\cdot$}}

\put(10,98){\vector(0,-1){25}}
\put(40,98){\vector(0,-1){25}}
\put(70,98){\vector(0,-1){25}}
\put(100,98){\vector(0,-1){25}}
\put(130,98){\vector(0,-1){25}}

\put(60,125){\vector(-2,-1){50}}
\put(65,125){\vector(-1,-1){25}}
\put(70,125){\vector(0,-1){25}}
\put(75,125){\vector(1,-1){25}}
\put(80,125){\vector(2,-1){50}}

\end{picture}
\end{center}

\noindent  are equal amounts to proving that the space between all
these paths could be filled in by a complex of commutative
diagrams homeomorphic to an $n$-dimensional sphere. Such is, for
example, the following complex, called the \emph{associahedron},
or \emph{Stasheff polytope}:

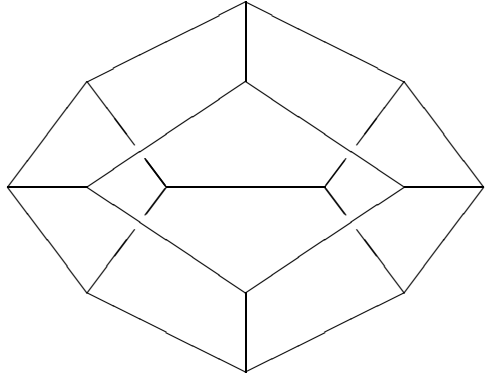
\begin{center}
\begin{picture}(180,160)

\put(90,10){\line(-2,1){60}}

\put(90,10){\line(2,1){60}}

\put(90,10){\line(0,1){30}}

\put(30,40){\line(-3,4){30}}

\put(30,40){\line(3,4){18}}

\put(150,40){\line(-3,4){18}}

\put(150,40){\line(3,4){30}}

\put(90,40){\line(-3,2){60}}

\put(90,40){\line(3,2){60}}

\put(0,80){\line(3,4){30}}

\put(0,80){\line(1,0){30}}

\put(180,80){\line(-3,4){30}}

\put(180,80){\line(-1,0){30}}

\put(60,80){\line(-3,4){8}}

\put(60,80){\line(-3,-4){8}}

\put(60,80){\line(-3,4){8}}

\put(60,80){\line(1,0){60}}

\put(120,80){\line(3,4){8}}

\put(120,80){\line(3,-4){8}}

\put(30,80){\line(3,2){60}}

\put(150,80){\line(-3,2){60}}

\put(30,120){\line(3,-4){18}}

\put(30,120){\line(2,1){60}}

\put(150,120){\line(-3,-4){18}}

\put(150,120){\line(-2,1){60}}

\put(90,120){\line(0,1){30}}

\end{picture}
\end{center}

\noindent whose vertices are all five-letter terms made with one
binary operation, and whose edges correspond to single
applications of the associativity law. Then the equality of two
paths follows from the fact that they are homotopic in the
complex. This is the \emph{global} approach to coherence, which
stems from \cite{St63} (see also Stasheff's papers in \cite{LO97},
and references therein).

There is also a \emph{local} approach to coherence, which stems
form \cite{ML63}. In the term-rewriting terminology, we have to
prove that for any two paths of arrows that terminate in the same
normal form:

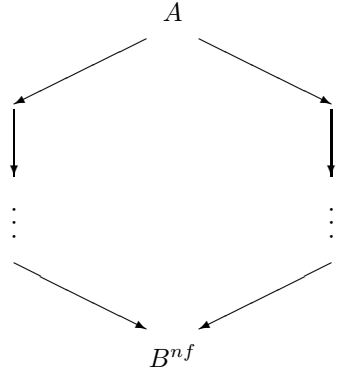
\begin{center}
\begin{picture}(140,140)
\put(70,5){\makebox(0,0){$B^{nf}$}}
\put(70,135){\makebox(0,0){$A$}}

\put(10,40){\vector(2,-1){50}}
\put(130,40){\vector(-2,-1){50}}

\put(10,50){\makebox(0,0){$\cdot$}}
\put(10,55){\makebox(0,0){$\cdot$}}
\put(10,60){\makebox(0,0){$\cdot$}}

\put(130,50){\makebox(0,0){$\cdot$}}
\put(130,55){\makebox(0,0){$\cdot$}}
\put(130,60){\makebox(0,0){$\cdot$}}

\put(10,98){\vector(0,-1){25}}
\put(130,98){\vector(0,-1){25}}

\put(60,125){\vector(-2,-1){50}}
\put(80,125){\vector(2,-1){50}}

\end{picture}
\end{center}

\noindent one can tile the space in between by commuting diagrams
of arrows (reductions). For this tiling we proceed inductively in
the following manner (see \cite{JohnM87}, Lemma 4.3, where the
assumption that we deal with isomorphisms is replaced by the
weaker assumption that we deal with monomorphisms; cf.\ also
\cite{DP04}, Section 4.3):

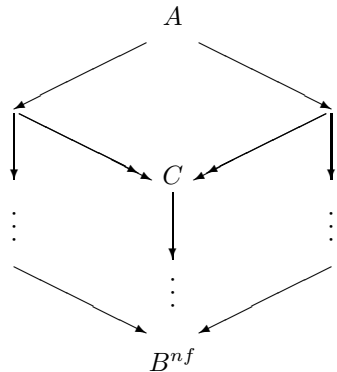
\begin{center}
\begin{picture}(140,140)
\put(70,5){\makebox(0,0){$B^{nf}$}}
\put(70,135){\makebox(0,0){$A$}}
\put(70,75){\makebox(0,0){$C$}}

\put(10,40){\vector(2,-1){50}}
\put(130,40){\vector(-2,-1){50}}

\put(10,50){\makebox(0,0){$\cdot$}}
\put(10,55){\makebox(0,0){$\cdot$}}
\put(10,60){\makebox(0,0){$\cdot$}}

\put(130,50){\makebox(0,0){$\cdot$}}
\put(130,55){\makebox(0,0){$\cdot$}}
\put(130,60){\makebox(0,0){$\cdot$}}

\put(70,25){\makebox(0,0){$\cdot$}}
\put(70,30){\makebox(0,0){$\cdot$}}
\put(70,35){\makebox(0,0){$\cdot$}}

\put(10,98){\vector(0,-1){25}}
\put(130,98){\vector(0,-1){25}}
\put(70,68){\vector(0,-1){25}}

\put(60,125){\vector(-2,-1){50}}
\put(80,125){\vector(2,-1){50}}

\put(128,98){\vector(-2,-1){50}}
\put(128,98){\vector(-2,-1){45}}
\put(12,98){\vector(2,-1){50}}
\put(12,98){\vector(2,-1){45}}

\end{picture}
\end{center}

\noindent At this place, we are faced with all the difficulties
that appear in proofs of the Church-Rosser property for a notion
of reduction, which consist in listing all the critical pairs of
reductions. The difference with what we have in term-rewriting
systems is that we must always verify that our tiles are commuting
diagrams of arrows. In term-rewriting systems we usually do not
deal with that (but cf.\ \cite{M03}, and references therein; the
procedure sketched above works when all the paths starting from
the same vertex are bounded in length).

It is not however true that all the interesting cases of coherence
where ``all diagrams commute'' involve only arrows that are
isomorphisms (see \cite{Lap72b}, \cite{JohnM87}, Lemma 4.2, and
\cite{DP04}, Section 4.2; remark that the four-dimensional
associahedron has 42 vertices). Consider, for example, arrows
whose type

\[
A\wedge (B\vee C) \vdash (A\wedge B)\vee C
\]

\noindent has something to do both with distributivity and
associativity, and which in \cite{DP04} are called
\emph{dissociativity} arrows (in the literature, the same
principle is also called weak or linear distribution; see
\cite{DP04}, Section 7.1, for references from category theory,
logic and universal algebra). These arrows need not be
isomorphisms, and they are of particular interest because they
underlie the cut principle in multiple-conclusion (plural) sequent
systems.

If we have such arrows, which are not isomorphisms, then the
global approach to coherence is not open any more, and we have to
take the local approach. When we want to show that two paths of
arrows are equal by closing the initial forking of arrows by a
commutative diagram, as in the following picture:

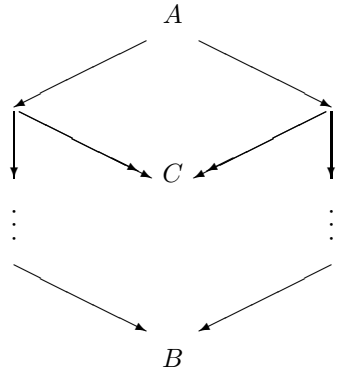
\begin{center}
\begin{picture}(140,140)
\put(70,5){\makebox(0,0){$B$}}
\put(70,135){\makebox(0,0){$A$}}
\put(70,75){\makebox(0,0){$C$}}

\put(10,40){\vector(2,-1){50}}
\put(130,40){\vector(-2,-1){50}}

\put(10,50){\makebox(0,0){$\cdot$}}
\put(10,55){\makebox(0,0){$\cdot$}}
\put(10,60){\makebox(0,0){$\cdot$}}

\put(130,50){\makebox(0,0){$\cdot$}}
\put(130,55){\makebox(0,0){$\cdot$}}
\put(130,60){\makebox(0,0){$\cdot$}}

\put(10,98){\vector(0,-1){25}}
\put(130,98){\vector(0,-1){25}}

\put(60,125){\vector(-2,-1){50}}
\put(80,125){\vector(2,-1){50}}

\put(128,98){\vector(-2,-1){50}}
\put(128,98){\vector(-2,-1){45}}
\put(12,98){\vector(2,-1){50}}
\put(12,98){\vector(2,-1){45}}

\end{picture}
\end{center}

\noindent we need an efficient criterion for showing that the
object $C$ is still ``above'' the object $B$ (here $B$ need not to
be in normal form); i.e., we need to show that we have a path of
arrows from $C$ to $B$:

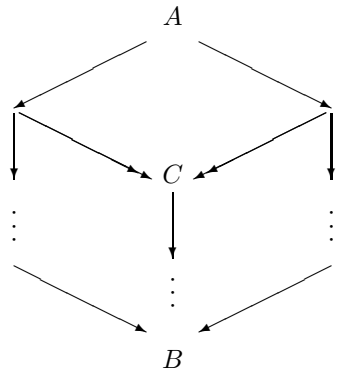
\begin{center}
\begin{picture}(140,140)
\put(70,5){\makebox(0,0){$B$}}
\put(70,135){\makebox(0,0){$A$}}
\put(70,75){\makebox(0,0){$C$}}

\put(10,40){\vector(2,-1){50}}
\put(130,40){\vector(-2,-1){50}}

\put(10,50){\makebox(0,0){$\cdot$}}
\put(10,55){\makebox(0,0){$\cdot$}}
\put(10,60){\makebox(0,0){$\cdot$}}

\put(130,50){\makebox(0,0){$\cdot$}}
\put(130,55){\makebox(0,0){$\cdot$}}
\put(130,60){\makebox(0,0){$\cdot$}}

\put(70,25){\makebox(0,0){$\cdot$}}
\put(70,30){\makebox(0,0){$\cdot$}}
\put(70,35){\makebox(0,0){$\cdot$}}

\put(10,98){\vector(0,-1){25}}
\put(130,98){\vector(0,-1){25}}
\put(70,68){\vector(0,-1){25}}

\put(60,125){\vector(-2,-1){50}}
\put(80,125){\vector(2,-1){50}}

\put(128,98){\vector(-2,-1){50}}
\put(128,98){\vector(-2,-1){45}}
\put(12,98){\vector(2,-1){50}}
\put(12,98){\vector(2,-1){45}}

\end{picture}
\end{center}

\noindent A criterion for the existence of such a path in the case
where we have associativity isomorphisms and dissociativity
arrows, which are not isomorphisms, is spelled out in \cite{DP04}
(Section 7.3, Theoremhood Proposition).

Coherence in this case could perhaps also be deduced from a very
general theorem of \cite{Ben68} (Theorem 5.2.4), whose proof is
only sketched in that paper, with substantial parts missing. It is
not clear whether the proof of \cite{DP04} (Section 7.3) was
envisaged in \cite{Ben68}, and judging by the complexity of
particular criteria, as the one mentioned in the preceding
paragraph, this seems unlikely.

In cases where such a criterion is not available, the paths of
arrows should first be normalized, according to some normalization
procedure (this is often a procedure inspired by cut elimination),
and then, in order to establish coherence, one has to compare such
normalized arrows (see, for example, \cite{KML71}, \cite{DP04},
Chapters 7-14, and \cite{DP05}). The normal form of paths of
arrows need not be unique.


\begin{thebibliography}{99}

\bibitem{Ben63} {\sc J. B\' enabou},
{\it Cat\' egories avec multiplication}, \textbf {\textit {Comptes
Rendus de l'Acad\' emie des Sciences, Paris}}, S\' erie I, Math\'
ematique, vol.\ 256 (1963), pp.\ 1887-1890

\bibitem{Ben68} --------,
{\it Structures alg\' ebriques dans les cat\' egories}, \textbf
{\textit {Cahiers de Topologie et G\' eom\' etrie Diff\'
erentielle}}, vol. 10 (1968), pp. 1-126

\bibitem{BZ85} {\sc G. Burde} and {\sc H. Zieschang}, \textbf {\textit {Knots}},
de Gruyter, Berlin, 1985

\bibitem{CK73} {\sc C.C. Chang} and {\sc H.J. Keisler}, \textbf {\textit {Model Theory}},
North-Holland, Amsterdam, 1973

\bibitem{DP04} {\sc K. Do\v sen} and {\sc Z. Petri\' c}, \textbf {\textit
{Proof-Theoretical Coherence}}, KCL Publications, London, 2004

\bibitem{DP05} --------, \textbf {\textit
{Proof-Net Categories}}, preprint, Mathematical Institute,
Belgrade, 2005

\bibitem{Huet80} {\sc G. Huet},
{\it Confluent reductions: Abstract properties and applications to
term rewriting systems}, \textbf {\textit {Journal of the
Association for Computing Machinery}}, vol.\ 27 (1980), pp.
797-821

\bibitem{JohnM87} {\sc M. Johnson}, \textbf {\textit {Pasting Diagrams in n-Categories
with Applications to Coherence Theorems and Categories of Paths}},
doctoral thesis, University of Sydney, 1987

\bibitem{JohnM92} --------, {\it Linear term rewriting systesm are higher
dimensional string
rewriting systems}, \textbf {\textit {Proceedings of the Institute
for Mathematics and its Applications}} (C.M.I. Rattray and R.G.
Clark, editors), vol.\ 35, Oxford University Press, 1992, pp.
101-110

\bibitem{KML71} {\sc G.M. Kelly} and {\sc S. Mac Lane},
{\it Coherence in closed categories}, \textbf {\textit {Journal of
Pure and Applied Algebra}}, vol.\ 1 (1971), pp.\ 97-140, 219

\bibitem{LS86} {\sc J. Lambek} and {\sc P.J. Scott},
\textbf {\textit {Introduction to Higher Order Categorical
Logic}}, Cambridge University Press, Cambridge, 1986

\bibitem{Lap72b}  {\sc M.L. Laplaza}, {\it Coherence for associtivity not an isomorphism},
\textbf {\textit {Journal of Pure and Applied Algebra}}, vol. 2
(1972), pp. 107-120

\bibitem{LO97} {\sc J.-L. Loday} et al., editors,
\textbf {\textit {Operads: Proceedings of Renaissance
Conferences}}, Contemporary Mathematics, vol. 202, American
Mathematical Society, Providence, 1997

\bibitem{ML63} {\sc S. Mac Lane}, {\it Natural associativity and
commutativity}, \textbf {\textit {Rice University Studies, Papers
in Mathematics}}, vol.\ 49 (1963), pp.\ 28-46

\bibitem{ML71} --------, \textbf {\textit {Categories for the Working
Mathematician}}, Springer, Berlin, 1971 (expanded second edition,
1998)

\bibitem{M03} {\sc P.-A. Melli\` es}, {\it Axiomatic rewriting theory VI:
Residual theory revisited}, \textbf {\textit {Rewriting Techniques
and Applications}} (S. Tison, editor), Lecture Notes in Computer
Science, vol. 2378, Springer, Berlin, 2002, pp.\ 24-50

\bibitem{PD} {\sc A. Preller} and {\sc P. Duroux}, {\it
Normalisation of the theory \textbf{T} of Cartesian closed
categories and conservativity of extensions \textbf{T}$\,${\rm
[$x$]} of \textbf{T}}, \textbf {\textit {Theoretical Informatics
and Applications}}, vol.\ 33 (1999), pp.\ 227-257

\bibitem{St63} {\sc J.D. Stasheff}, {\it Homotopy associativity of H-spaces, I, II},
\textbf {\textit {Transactions of the American Mathematical
Society}}, vol. 108 (1963), pp. 275-292, 293-312

\bibitem{S75} {\sc M.E. Szabo}, {\it A counter-example to coherence
in cartesian closed categories}, \textbf {\textit {Canadian
Mathematical Bulletin}}, vol. 18 (1975), pp. 111-114


\end{thebibliography}
\end{document}